\documentclass[letterpaper,12pt,twoside]{article}
\usepackage[left=2.5cm, right=2.5cm,tmargin = 2.5cm]{geometry}
\usepackage{amsmath,amssymb,amsfonts}
\usepackage{indentfirst}
\usepackage[utf8]{inputenc}
\usepackage{hyperref}
\usepackage{tikz}
\usepackage{mathtools,amsthm}

\newtheorem{defi}{Definition}[section]
\newtheorem{teo}[defi]{Theorem}
\newtheorem{lem}[defi]{Lemma}
\newtheorem{cor}[defi]{Corollary}
\newtheorem{prop}[defi]{Proposition}

\def\R{\mathbb R}

\usetikzlibrary{shapes.geometric,shapes,arrows,decorations,optics,graphs,babel,positioning}

\newcommand{\pathsimp}
{
\begin{tikzpicture}[node distance={32mm},vrt/.style = {draw,circle,minimum size=0.85cm}] 
\node[vrt] (1a) {$v_1$};
\node[vrt] (2a) [right of = 1a]{};
\node[vrt] (3a) [right of = 2a]{};
\node (pts) [right of = 3a, node distance = 15mm]{$\cdots$};
\node[vrt] (1b) [right of = pts, node distance = 15mm]{};
\node[vrt] (2b) [right of = 1b]{};
\node[vrt] (3b) [right of = 2b]{$\small v_n$};
\draw (1a) -- (2a);
\draw (2a) -- (3a);
\draw (1b) -- (2b);
\draw (2b) -- (3b);
 \node (ce1) [right of = 1a, node distance = 13mm]{};
\node (ce2) [right of = 2a, node distance = 13mm]{};
\node (ce3) [right of = 1b, node distance = 13mm]{};
\node (ce4) [right of = 2b, node distance = 13mm]{};
\node (c1) [above of = ce1, node distance = 3mm]{$c_1$};
\node (c2) [above of = ce2, node distance = 3mm]{$c_2$};
\node (c3) [above of = ce3, node distance = 3mm]{$c_{n-2}$};
\node (c4) [above of = ce4, node distance = 3mm]{$c_{n-1}$};
\end{tikzpicture}
}

\title{Maximal Algebraic Connectivity for Paths with Fixed Effective Resistance}
\author{Alonso Cruz Ortega \and Federico Menéndez--Conde}
\begin{document}
\maketitle

\begin{abstract}
We consider the problem on finding the edge weights that maximize the algebraic connectivity of a graph, subject to the condition that the total effective resistance is kept constant. We propose the conjecture that for every graph the maximum is attained for weights that are invariant under automprphisms. The solution to the problem is given explicitly for the paths $P_3$ y $P_4$, where the conjecture holds. 

\end{abstract}

\section{Introduction}
\medskip 

In 1973 M. Fiedler used the expression {\it algebraic connectivity} for the second smallest eigenvalue of the Laplacian of a graph \cite{fiedler73}. Notorious relations between that eigenvalue and the geometric idea of connectivity of the graph have been shown to hold ever since. Fiedler called the corresponding eigenvectors by the name {\it characteristic valuations} of the graph. There is an extensive use by many authors of the names {\it Fiedler value} for the algebraic connectivity and {\it Fiedler vectors} for the characteristic valuations. We refer to \cite{abreu, mohar} for historical overviews and many results on this matter. On top of their theoretical interest, Fiedler values and vectors have been studied in numerous applications, such as spectral clustering 
\cite{b-k,DP-S}, parallel processing \cite{simon}, biological evolution \cite{n-h}, neuroscience \cite{b-p}, air transportation \cite{wei}, among others. 
In general, the algebraic connectivity might be of interest whenever one would like to understand the global connectedness of a network.
\medskip 

In \cite{fiedler75} the definition of algebraic connectivity was extended to include graphs with edge weights, and in \cite{fiedler90} Fiedler introduced the {\it absolute algebraic connectivity}, defined as the maximum possible value of the algebraic connectivity of a weighted graph under the assumption that the sum of the edge weights equals the number of edges. In the same reference, Fiedler found the absolute algebraic connectivity of trees and showed that for every graph the absolute algebraic connectivity is attained for weights that are invariant under automorphisms. As a consequence of this, the absolute connectivity of  cycles equals tha algebraic connectivity of the cycle with all weights equal to 1 i.e. the {\it unweighted cycle}.
 \medskip 

The trace of the Laplacian of a graph equals the sum of the degrees of the vertices, hence twice the sum of the weights of the edges. Therefore, finding the absolute algebraic connectivity is equivalent to the problem on maximizing the algebraic connectivity of the graph under the condition that the trace of the Laplacian is constant, that is:   
\[
\lambda_2+ \cdots +\lambda_n = 2|E|,
\]
where $\lambda_j$ are the ordered eigenvalues of the Laplacian and $|E|$ is the number of edges. We note that the right-hand side in the equality above corresponds to the trace of the unweighted graph. In this work we consider the analogous problem of 
determining the weights that maximize the algebraic connectivity of a graph subject to the condition that the sum
\[
\frac 1 {\lambda_2(G;\omega)} + \cdots + \frac 1 {\lambda_n(G;\omega)}
\]
is held constant. This is known as the {\it Kirchhoff index}, and equals the trace of the Moore - Penrose pseudo-inverse of the Laplacian. In \cite{k-r}, Klein and Randic showed that this sum is related to the {\it total effective resistance} of the graph $R(G)$ by the identity
\begin{equation}\label{KI}
R(G) = n \sum_{j=2}^n \frac 1 {\lambda_n}, 
\end{equation}

; this is defined as the sum 
of the distances of all pairs of vertices of the graph with respect to the {\it effective resistance metric}. This metric decreases when weights increase or edges are added to the graph. This property makes this metric a natural choice in situation where more connections should mean more closedness.  The study of this metric originated in electric circuit theory \cite{d-s} and, in view of its aforementioned properties,  it has been considered in diverse applications such as molecular structures \cite{bklnt}, air traffic \cite{ymqw}, discrete geometry \cite{d-l}, network clustering \cite{vnls}, among others. This metric has also been defined for some infinite sets, and plays an important role in the constructions of Laplacians in self-similar fractals
\cite{kigami,strichartz} and in infinite neetworks \cite{j-p}. 

In this work we are interested in the relation between these two quantities: the algebraic connectivity and the total effective resistance of the graph. As we mentioned above, both indices are related to the connectedness of the network. In particular, we are interested in finding the edge weights that maximize the algebraic connectivity under the condition that the total effective resistance equals a given constant.
By the identity \eqref{KI}, the problem we consider is that of maximizing $\lambda_2$ while the trace of the pseudo-inverse of the Laplacian is fixed. In this sense, the proposed problem can be viewed as an analogue of the problem of finding the absolute algebraic connectivity.
It is therefore natural to consider the question if for our problem the solution will also be attained for weights that are 
invariant under automorphisms. We believe that this is likely true in the general case. In this work we solve the problem explicitly for the paths $P_3$ and $P_4$. Unfortunately our methods do not generalize in an obvious way to general $P_n$, let alone to more general graphs. 
The results are included in section \ref{principal}.
\bigskip

\section{Notation and Preliminaries}\label{preliminares}
\medskip 

\begin{figure}[h]
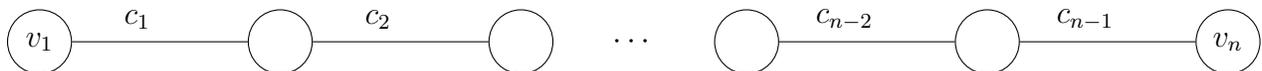

\begin{center}
\pathsimp
\end{center}
\caption{The path $P_n$ with weights.}\label{path}
\end{figure}

A {\it graph} $G$ is a set of vertices $V=(v_1,\dots,v_n)$ which we will always consider ordered, together with a set of unordered pairs of vertices (the {\it edges} of the graph). The set of edges is denoted by $E$. We assign weights to the edges by means of a function $\omega:E\rightarrow(0,\infty)$ and
write $\omega_{i,j}=\omega(\{v_i,v_j\})$. The pair $(G;\omega)$ is a {\it weighted graph}. The {\it unweighted graph} corresponds to the case when $\omega_{i,j}=1$ at every edge. We write $P_n$ for the {\it path} with $n$ vertices: that is the graph with  $V=(v_1,\dots, v_n)$ end $E=\{\{v_i,v_{i+1}\}\}_{i=1}^{n-1}$. We underline that the vertices in $P_n$ are ordered from one end-point to the other by neighbouring vertices. The weights of $P_n$ are denoted by $c_i=\omega_{i,i+1}$, as in figure \ref{path}. We will also use the notation $r_i=1/\omega_i$. Here $r_i$ is the resistance distance between consecutive vertices of the graph $v_i$ y $v_{i+1}$ (see section \ref{mre}). 
\medskip 

\medskip 

We denote by $\ell(V)$ the space of real-valued functions with domain $V$. For $x\in\ell(V)$ we write $x=(x_1,\dots, x_n)$ where $x_j=x(v_j)$. In the vector space $\ell(V)$ we consider the (real) inner product $\bigl<x,y\bigr>=x_1 y_1 +\cdots+ x_n y_n$ and write $x\perp y$ to denote orthogonallity with respect to this product.
\medskip 

The {\it laplacian} of a weighted graph $(G;\omega)$ is the operator on $\ell (V)$ given by 
\[
(L_{G,\omega}\, x)_i = \sum_{\{v_i,v_j\}\in E} \omega_{i,j}(x_i-x_j).
\]
An equivalent way to define the Laplacian is as the self-adjoint operator determined by the quadratic form
\[
\bigl<L_{G,\omega}\, x,x\bigr>=\sum_{\{v_i,v_j\}\in E} \omega_{i,j} (x_i-x_j)^2. 
\]
When it causes no confusion, we will simply write $L_\omega$ instead of $L_{G,\omega}$. If $\omega$ is not in the subindex, we are considering the unweighted graph. 
\medskip 

The graph Laplacian is a singular non-negative operator, and the constant vector $J=(1,1,\dots,1)$ is always in its kernel. 
The eigenvalues of $L_{G,\omega}$ are written in non-decreasing order as $0=\lambda_1(G;\omega)\leq\lambda_2(G;\omega)\leq\cdots\leq\lambda_n(G;\omega)$.  
\medskip 

The Moore--Penrose pseudo-inverse of $L_\omega$, denoted $L_\omega^+$, is determined by the condition that $L_\omega^+ L_\omega = L_\omega L_\omega^+$ is the orthogonal projection over the subspace $({\rm Ker\  }L_\omega)^\perp$, which for connected graphs is given by $J^\perp$, a subspace of codimension 1.
\bigskip

\section{The algebraic connectivity}\label{con_alg}
\medskip 

By the Courant-Fischer-Weyl theorem, the algebraic connectivity of a weighted graph is 
\[
\lambda_2(G;\omega) = \min_{x\perp J} 
\frac {\bigl<L_{G,\omega}\, x,x\bigr>} {\bigl<x,x\bigr>}.
\]
The quotient in the right hand side is the {\it Rayleigh's quotient} of the operator $L_{G,\omega}$ valued at $x$. A well known consequence of this 
is that $\lambda_2>0$ if and only if the graph is connected. Also, $\lambda_2$ cannot increase if weights are increased or if edges are added.
\medskip 

Nodal domains for eigenvectors of the graph Laplacians have been studied extensively  (e.g. \cite{fiedler75,f-k}).
For the particular case of Fiedler vectors in paths, the following condition is satisfied:
\medskip 

\begin{prop}\label{monotono}
If $u$ is a Fiedler vector of $(P_n;\omega)$, then $u$ is either increasing or decreasing. 
\end{prop}
\medskip 

Even more can be said when the weights of $P_n$ are symmetric. The next proposition characterizes the Fiedler vectorsin this case with $n$-even. 
\medskip 

\begin{prop}\label{simetria}
Suppose that in $(P_n;\omega)$ the weights $\{c_1,\dots,c_{n-1}\}$ are symmetric in the sense that $c_i=c_{n-i}$ for all $i=1,\dots,n-1$. Then, for a Fiedler vector 
$x=(x_1,\dots,x_n)$ the following holds: 
\[
x_j+x_{n+1-j}=0,\qquad\qquad j=1,\dots, n.
\]
\end{prop}

\begin{proof}
Let $L_\omega$ be the Laplacian for $(P_n;\omega)$, and $x$ a corresponding Fiedler vector.
Introduce a new vector $x'=(x_1',\dots,x_n')$, where $x_j'=x_{n+1-j}$. 
Making $c_0=c_n=x_0=x_{n+1}=0$ for convenience, we can write for all $j=1, 2, \dots, n$ 
\begin{align*}
(L_\omega\, x')_j & = c_{j-1}(x_j'-x_{j-1}')+c_j(x_j'-x_{j+1}')\cr
& = c_{n-j+1}(x_{n-j+1}-x_{n-j+2})+c_{n-j}(x_{n-j+1}-x_{n-1})\cr
& = (Lx)_{n-j+1}\cr
& = \lambda_2\hskip0.05cm x_{n+1-j}\cr
& = \lambda_2\hskip0.05cm x_j'.
\end{align*}
So,  $x'$ is an eigenvector for $\lambda_2$. Since the Laplacian matrix is tridiagonal, symmetric and the entries next to the main diagonal are different from zero, it can be whown that all of its eigenvalues are simple. It follows that  $x'$ and $x$ are linearly dependent, say  $x'=\alpha x$. This implies that $x_{n+1-j}=\alpha x_j$, and since the entries of $x$ cannot have constant sign, necessarily $\alpha<0$. Then,

\begin{align*}
0 &= x_1 +\cdots +x_n\cr
&=(1+\alpha)\sum_{j=1}^{\frac n 2} x_j.
\end{align*}

The $x_j$ in the sum above have the same sign, so that we have $1+\alpha=0$ which gives the result. 
\end{proof}
\bigskip 

\section{The effective resistance metric}\label{mre}
\medskip 

Let  $(G;\omega)$ be a connected weighted graph and $U\subset V$ a subset of its vertices. Given $x\in\ell(U)$ there exists a unique extension of $x$ to $V$, which we denote by $\tilde x\in\ell(V)$, such that  $\left(L_\omega\, \tilde x\right) v=0$ for all  $v\in V\setminus U$. The vector $\tilde x$ is known as the {\it harmonic 
extension} of $x$. 
\medskip

\begin{defi}\label{mref}
Given two vertices $v_i\neq v_j$ of the weighted graph $(G;\omega)$, let $U=\{v_i,v_j\}$. For $h\in\ell (V)$ the harmonic extension of the function in $\ell(U)$ given by $h_i=1$, $h_j=0$, define
\[
r(v_i,v_j) = \bigl<L_\omega h,h\bigr>^{-1} = (L_\omega h)_i^{-1}.
\]
The number $r(v_i,v_j)$ is the {\it effective resistance} between $v_i$ and $v_j$. Also define  $r(v,v)=0$ for every vertex $v$. 
\end{defi}

It can be shown that $r(\cdot,\cdot)$ defines a metric on the set of vertices (e.g. \cite{kigami}). Different equivalent formulations for the resistance metric are known. One of those is given by the quadratic form 

\begin{equation}\label{psin}
r(v_i,v_j) = \bigl<L^+_\omega (\chi_i-\chi_j),\chi_i-\chi_j\bigr>
\end{equation}
where $L^+$ is the Moore-Penrose pseudo-inverse of the Laplacian, and $\chi_k$ is the characteristic function of $\{v_k\}$.
\medskip 

The total effective resistance is just the sum of all the distances determined by the effective resistance metric.
\smallskip 

\begin{defi}
The total effective resistance of $(G,\omega)$ is
\[
R(G;\omega)=\sum_{\{v_i,v_j\}} r(v_i,v_j).
\]
\end{defi}

We underline that the sum in the definition above involves the distances between every pair of vertices, regardless on whether they share an edge. 
\medskip 

Using \eqref{psin}, in \cite{k-r} it was shown that
\[
R(G;\omega)=n\sum_{j=2}^n \frac 1 {\lambda_j(G;\omega)}.
\]

Note that the right hand side equals the trace of the pseudo-inverse $L_\omega^+$.
\medskip 

The following result is well known (e.g. \cite{d-s}) and it will be useful to find the total effective resistance of paths $P_n$. We include the proof for 
self-reference.
\medskip 

\begin{teo}\label{enserie}
Let $(G,\omega)$ be a connected weighted graph such that there exist two subgraphs $\Gamma_1$ and $\Gamma_2$, that share a single vertex $v$. Suppose that there are no edges joning vertices in $\Gamma_1\setminus\{v\}$ with vertices in $\Gamma_2\setminus\{v\}$. Let $w_1$ and $w_2$ be vertices in $G_1$ and $G_2$, respectively. Then  
\[
r(w_1,w_2) = \rho_1+\rho_2,
\]
where $\rho_j$ is the effective resistance in $\Gamma_j$ between $w_j$ and $v$.
\end{teo}

\begin{proof}
For $j=1,2$, let $\varphi_j$ be harmonic in $G\setminus \{w_j,v\}$,
with $\varphi_j(w_j)=1$ and $\varphi_j(v)=0$. Note that $\varphi_1=0$ in $\Gamma_2$ and $\varphi_2=0$ in $\Gamma_1$. From definition\ref{mref}, we have
\[
\rho_j=(L\varphi_j)(w_j)^{-1}.
\]

We want to find a function $\psi$, harmonic outside the set $\{w_1,w_2\}$, of the form 
\[
\psi = \alpha_1\varphi_1+\alpha_2\varphi_2,
\]
with $\alpha_j$ some constants. Noting 
\begin{align*}
(L\psi)(v) & = \alpha_1 (L\varphi_1)(v)+\alpha_2 (L\varphi_2)(v)\cr
& = -\frac{\alpha_1}{\rho_1}-\frac{\alpha_2}{\rho_2}
\end{align*}
we can take
\[
\alpha_1 = -\frac 1 {\rho_2},\qquad\qquad \alpha_2 = \frac 1 {\rho_1}.
\]
\medskip 

From $\psi(w_j)=\alpha_j$, it follows that $\varphi$ definided by
\[
\varphi = \frac{\psi-\alpha_2}{\alpha_1-\alpha_2}
\]
is the harmonic extension to $G$ of the function in  $\{w_1,w_2\}$ given by $\varphi(w_1)=1$ and $\varphi(w_2)=0$.
\medskip 

We obtain
\begin{align*}
(L\varphi)(w_1) & = \frac{(L\psi)(w_1)}{\alpha_1-\alpha_2}\cr
& = \frac{\alpha_1 (L\varphi_1)(w_1)+\alpha_2 (L\varphi_2)(w_1)}{\alpha_1-\alpha_2}\cr
& = \frac{\alpha_1 (L\varphi_1)(w_1)}{\alpha_1-\alpha_2}.
\end{align*}
By a simple substitution it follows 
\[
\frac 1 {r(w_1,w_2)} = \frac{\frac{1}{\rho_1}\frac{1}{\rho_2}}{\frac{1}{\rho_1}+\frac{1}{\rho_2}}
\]
from which we can see that $r(w_1,w_2)=\rho_1+\rho_2$.
\end{proof}
\medskip 

The following corollary, stating that the resistances between vertices in a path add up is an immediate consequence of theorem \ref{enserie}. 
\smallskip 

\begin{cor}
Let $(P_n;\omega)$ be a path with weights $\{c_1,\dots, c_{n-1}\}$ and let $v_i, v_j$ vertices in $P_n$, with $i<j$. The effective resistance between $v_i$ and $v_j$ is given by
\[
r(v_i,v_j)  = \sum_{k=i}^{j-1} \frac 1 {c_k} = \sum_{k=i}^{j-1} r(v_k,v_{k+1}).
\]
\end{cor}
\bigskip 

\section{Algebraic connectivity in paths}\label{principal}
\medskip 

Define 

\begin{equation}\label{breve}
\breve a (G) = \max\left\{\lambda_2(G;\omega)\ |\ R(G;\omega)=R(G)\right\}.
\end{equation}
\medskip 

We conjecture that -- in the same fashion of what happens with the absolute algebraic connectivity -- the maximal value of \eqref{breve} is attained always for weight distributions that are invariant under automorphisms of the graph. In what follows, we will show that this is the case for the paths $P_3$ and $P_4$, finding explicitly the maximum in both cases.
\medskip 

From the definition, it follows that $r(v_i,v_j)$ is inversely proportional to multiplcation of the weights by a constant. Therefore
\[
R(G;c\hskip0.08cm\omega) = \frac 1 c R(G;\omega),
\]
which can also be seen from the identity \eqref {KI}. So, for every weighted graph $(G;\omega)$ the product $\lambda(G;c\, \omega) R(G;c\, \omega)$ is independent of $c$. The problem of finding $\breve a$ is equivelent to determining the weights $\omega$ that maximize the product $\lambda(G;\omega) R(G;\omega)$. 
\medskip 

Next, we show the solution for the path $P_3$. 
\medskip 

\begin{figure}
\begin{center}
\includegraphics[scale=0.59]{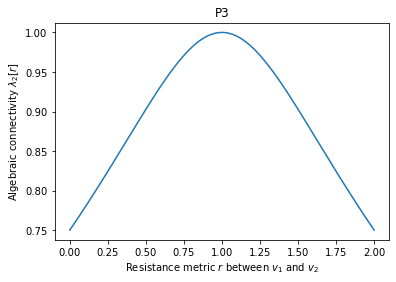}
\caption{Algebraic connectivity for $P_3$ with total effective resistance $r_1+r_2=2$.}
\label{fig:coneP3}
\end{center}
\end{figure}

\begin{teo}\label{p3}
Consider $(P_3;\omega)$ with weights $\{c_1,c_2\}$ satisfying
\[
R(P_3;\omega)=\frac 1 {c_1} + \frac 1 {c_2} = 2. 
\]
Then 
\[
\lambda_2(P_3;\omega)\leq \lambda_2(P_3),
\]
with equality iff $c_1=c_2=1$. That is
\[
\breve a (P_3) = \lambda_2(P_3) = 1.
\]
\end{teo}
\begin{proof} 
Let $u=(1+t,1-t,-2)$ be an eigenvector of $L_{\omega}$ with eigenvalue $\lambda$. From the equalities $(L_\omega u)_1 = 2c_1 t$ 
and $(L_\omega u)_3 = c_2 (t-3)$ it follows that
\begin{equation}\label{eclambda}
\lambda = \frac{2 c_1 t}{1+t} = \frac {c_2(3-t)} {2}.
\end{equation}

using the notation $b=c_2/c_1$, the second equality can be written as
\[
bt^2+(4-2b)t-3b,
\]
with roots 
\begin{equation}\label{lastes}
t_\pm =\frac{b-2\pm 2\sqrt{b^2-b+1}} {b}.  
\end{equation}
Substitution in \eqref{eclambda} shows that the two positive eigenvalues of $L_\omega$ are given by
\begin{equation}\label{lambdamaso}
\lambda_2 = \frac{c_2(3-t_\pm)}{2},\qquad \lambda_3 = \frac{c_2(3+t_\pm)}{2}.
\end{equation}
\medskip 

Let $r>0$ be the effective resistance between $v_1$ and $v_2$, i.e. $r=r_1=1/c_1$ and also $b=r/(2-r)$. Since $r$ determines the weights of the graph, we
have that the algebraic connectivity $\lambda_2$ is a function of $r$. We want to whow that $\lambda_2[r]$ attains its maximal value for  $r=1$. Substituting the expression for the sign $+$ of \eqref{lastes} in \eqref{lambdamaso}, after simplification we can see that 
\begin{align*}
\lambda_2[r] &= \frac {2-\sqrt{3r^2-6r+4}}{r(2-r)}\cr
&\ \cr
&=  \frac {{2-\sqrt{3r^2-6r+4}}}{r(2-r)}\cdot \frac{2+\sqrt{3r^2-6r+4}}
{2+\sqrt{3r^2-6r+4}}\cr
&\ \cr
&=  \frac{3}{2+\sqrt{3r^2-6r+4}}.
\end{align*}
Since the value $3r^2-6r+4$ is minimal for $r=1$, we conclude that $\lambda_2[r]$ attains its only maximal value also at $r=1$ (figure \ref{fig:coneP3}).
\end{proof}

{\sc Note:} the largest eigenvalue for the Laplacian on  $(P_3;\omega)$ corresponds to the choice of root $t_-$ en \eqref{lastes}, from which 
\[
\lambda_3[r] =  \frac{3}{2-\sqrt{3r^2-6r+4}}.
\]
We can verify the identity \ref{KI} in this case:
\[
\frac 1 {\lambda_2[r]}+\frac 1 {\lambda_3[r]} = \frac 4 3.
\]
\medskip 

The following result is immediate from \ref{p3}.
\medskip 

\begin{cor}
For every choice of weights $\omega$ in $P_3$ we have that
\[
R(P_3;\omega)\lambda_2(P_3;\omega)\leq 4,
\]
with equality iff the weights are constant.
\end{cor}
\medskip 

In order to find $\breve a(P_4)$ we will first show that the weights that gives the optimal value of $\lambda_2$ are symmetrically distributed. To this aim, we show that the algebraic connectivity increases when the weights are substituted by its harmonic mean. Note that the same thing happened in theorem \ref{p3}.
\medskip 

\begin{lem}\label{simetriaP4}
Let $\{c_1,c_2,c_3\}$ be the weights in the edges of $(P_4;\omega)$ and let $d$ be the harmonic mean of $c_1$ and $c_3$. If $\{d,c_2,d\}$ are the weights in $(P_4;\omega')$, then 
\[
\lambda_2(P_4;\omega)\leq \lambda_2(P_4;\omega').
\]
\begin{proof}
In view of proposition \ref{simetria}, we know that there is a Fiedler vector $u$ for $(P_4;\omega')$ of the form $u=(-1,-a,a,1)$, with $a\in(0,1)$. The Rayleigh quotient gives 
\begin{equation}\label{casim}
\lambda_2(P_4;\omega')=\frac{2d(a-1)^2+4c_2a^2}{2(1+a^2)}.
\end{equation}
We will find a vector in $J^\perp$ such that its Rayleigh quotient for $(P_4;\omega)$ is less than $\lambda_2(P_4;\omega')$, which 
will give the result as this is an upper bound for $\lambda_2(P_4;\omega)$. 
\medskip 

Set
\[
v[\alpha] = (-1,-a,a,1) + \alpha(1,-1,-1,1),
\]
so that 
\begin{equation}\label{rayalfa}
\frac{\bigl<L_\omega v[\alpha],v[\alpha]\bigr>}{\|v[\alpha]\|}=
\frac{4\alpha^2(c_1+c_3)+4\alpha (1-a)(c_3-c_1)+(1-a)^2(c_1+c_3)+4c_2 a^2}{2(1+a^2)+4\alpha^2}.
\end{equation}
Since the denominator in \eqref{rayalfa} is less than the denominator in \eqref{casim}, it is enough to find an $\alpha$ for which the numerator in  \eqref{rayalfa} is not greater than the one in \eqref{casim}. We set $\alpha=\beta(1-a)$, so that the inequality that we want to verify can be written more concisely as
\begin{equation}\label{comparays}
(1+4\beta^2)(c_1+c_3)+4\beta(c_3-c_1)\leq 2d.
\end{equation}
The harmonic mean can be written as
\begin{align*}
d&=\frac{2}{\frac{1}{c_1}+\frac{1}{c_3}}=\frac{2c_1 c_3}{c_1+c_3}\cr
&=\frac{(c_1+c_3)^2-(c_1-c_3)^2}{2(c_1+c_3)}\cr
&=\frac{c_1+c_3}{2}-\frac{(c_1-c_3)^2}{2(c_1+c_3)}.
\end{align*}
Substitution in \eqref{comparays} gives
\[
4\beta^2(c_1+c_3)+4\beta(c_3-c_1)\leq -\frac{(c_1-c_3)^2}{c_1+c_3},
\]
which is solvable for $\beta\in\R$. For instance, taking
\[
\beta=\frac 1 2\, \frac{c_1-c_3}{c_1+c_3}
\]
gives the equality. 
\end{proof}
\end{lem}
\medskip 

\begin{cor}
If $\omega=\{c_1,c_2,c_3\}$ are weights for $(P_4;\omega)$ such that $\breve a(P_4) = \lambda_2(P4;\omega)$, then $c_1=c_3$.
\end{cor}
\begin{proof}
It follows from \ref{simetriaP4}, noting that 
 $(P_4;\omega)$ and  $(P_4;\omega')$ have the same total effective resistance.
\end{proof}
\medskip 

Now we use lemma \ref{simetriaP4} to find $\breve a (P_4)$. The maximum is attained for weights $\{d,c,d\}$ satisfying 
\begin{equation}\label{condi}
\frac 3 d + \frac 2 c =5.
\end{equation}
 
Taking the Fiedler vector $(-1,-t,t,1)$ $\lambda_2=\breve a (P_4)$, it follows that 
\begin{equation*}
\lambda_2=d(1-t)= d+2c-\frac d t.
\end{equation*}
Defining $b=(d+2c)/d$ in order to simplify the calculations, we obtain  
\[
t^2+(b+1)t-1=0.
\]
The root that corresponds to the lesser eigenvalue $\lambda_2$ is
\[
t = \frac {1-b} 2 + \frac{\sqrt{(b-1)^2+4}} 2,
\]
since we know that $0<t<1$ and $b>1$. Hence, the algebraic connectivity is in this case given by
\[
\lambda_2 = d\left(\frac{b+1} 2 - \frac{\sqrt{(b-1)^2+4}} 2\right).
\]
As we did for $P_3$, let $r$ de the effective resistance between $v_1$ and $v_2$. From \eqref{condi} we obtain the relations $d=1/r$, $c=2/(5-3r)$ and $b=(r+5)/(5-3r)$. The algebraic connectivity can be written in terms of $r\in(0,5/3)$. After calculations, we obtain: 
\[
\lambda_2[r] = \frac{1}{r(5-3r)}\left(5-r-\sqrt{13r^2-30r+25}\right).
\]

\begin{figure}
\begin{center}
\includegraphics[scale=0.59]{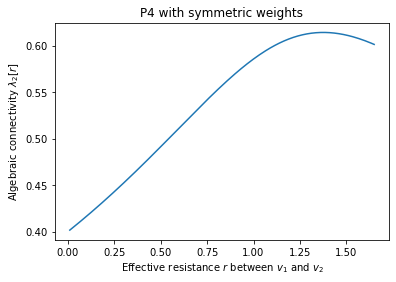}
\caption{Algebraic connectivity for $P_4$ with total effective resistance $10$ and $r_1=r_3$.}
\label{fig:coneP4}
\end{center}
\end{figure}

The value of $\lambda_2[r]$ is shown in figure \ref{fig:coneP4} where we can see that the maximum is attained at a single value of $r$. It is given by
\[
\lambda_2[r_0] = 2 - \frac {4\sqrt 3}{5},\qquad\qquad 
r_0=\frac{5}{13}\left(3+\frac 1 {\sqrt 3}\right).
\]

We have shown:
\medskip 

\begin{teo}
Let $(P_4;\omega_0)$ where tha weights $\{d_1,d_2,d_3\}$ of $\omega_0$ are given by
\[
d_1 = d_3 = \frac{13\sqrt 3}{5(\sqrt 3 -1)},\qquad 
d_2 = \frac5 {26}\left(6-(\sqrt 3 -1)\right). 
\]
If $(P_4;\omega)$ has weights $\{c_1,c_2,c_3\}$ such that
\[
\frac 3 {c_1} + \frac 4 {c_2} + \frac 3 {c_3} = 10,
\]
then 
\[
\lambda_2(P_4;\omega)\leq \lambda_2(P_4;\omega_0)= 2 - \frac {4\sqrt 3}{5},
\]
with equality iff $\omega=\omega_0$. In  particular 
$\breve a(P_4)=2 - \frac {4\sqrt 3}{5}$.
\end{teo}
 \medskip 
 
 \begin{cor}
The maximal value of the product $\lambda_2(P_4;\omega)R(P_4;\omega)$ is given by 
 $10(2 - \frac {4\sqrt 3}{5})$. 
 \end{cor}

Note that $\breve a(P_4)$ is irrational. This is different of what happens for the absolute algebraic connectivity, which is known to be rational for every 
tree \cite{fiedler90}.
\newpage

\newpage

\bibliography{referencias}{}
\bibliographystyle{ieeetr}
\label{LastPage}

\end{document}